\documentclass[12pt,reqno]{amsart}

\usepackage{amsmath}
\usepackage{amscd}
\usepackage{amssymb}
\usepackage{latexsym}
\usepackage{MnSymbol}
\usepackage{undertilde}
\usepackage{amsrefs}
\usepackage{nicefrac}

\input xy
\xyoption{all} 
\newdir{ >}{{}*!/-5pt/@{>}}

\theoremstyle{plain}
\newtheorem{theorem}{Theorem}[section]
\newtheorem*{theorem*}{Theorem}
\newtheorem{lemma}[theorem]{Lemma}
\newtheorem*{lemma*}{Lemma}
\newtheorem{corollary}[theorem]{Corollary}
\newtheorem{proposition}[theorem]{Proposition}
\newtheorem{defprop}[theorem]{Proposition-Definition}

\theoremstyle{remark}
\newtheorem{remark}[theorem]{Remark}

\theoremstyle{definition}
\newtheorem{definition}[theorem]{Definition}
\newtheorem{example}[theorem]{Example}

\makeatletter
\def\revddots{\mathinner{\mkern1mu\raise\p@
\vbox{\kern7\p@\hbox{.}}\mkern2mu
\raise4\p@\hbox{.}\mkern2mu\raise7\p@\hbox{.}\mkern1mu}}
\makeatother 
\newcommand{\bgl}{\begin{equation}} 
\newcommand{\egl}{\end{equation}}
\newcommand{\bgloz}{\begin{equation*}} 
\newcommand{\egloz}{\end{equation*}}
\newcommand{\bgln}{\begin{eqnarray}} 
\newcommand{\egln}{\end{eqnarray}}
\newcommand{\bglnoz}{\begin{eqnarray*}} 
\newcommand{\eglnoz}{\end{eqnarray*}}
\newcommand{\btheo}{\begin{theorem}}
\newcommand{\etheo}{\end{theorem}}
\newcommand{\btheooz}{\begin{theorem*}}
\newcommand{\etheooz}{\end{theorem*}}
\newcommand{\blemma}{\begin{lemma}}
\newcommand{\elemma}{\end{lemma}}
\newcommand{\blemmaoz}{\begin{lemma*}}
\newcommand{\elemmaoz}{\end{lemma*}}
\newcommand{\bproof}{\begin{proof}}
\newcommand{\eproof}{\end{proof}}
\newcommand{\bbew}{\begin{beweis}}
\newcommand{\ebew}{\end{beweis}}
\newcommand{\bremark}{\begin{remark}}
\newcommand{\eremark}{\end{remark}}
\newcommand{\bex}{\begin{example}\em}
\newcommand{\eex}{\end{example}}
\newcommand{\bdefin}{\begin{definition}}
\newcommand{\edefin}{\end{definition}}
\newcommand{\bprop}{\begin{proposition}}
\newcommand{\eprop}{\end{proposition}}
\newcommand{\bdefprop}{\begin{defprop}}
\newcommand{\edefprop}{\end{defprop}}
\newcommand{\bcor}{\begin{corollary}}
\newcommand{\ecor}{\end{corollary}}
\newcommand{\bfa}{\begin{cases}} 
\newcommand{\efa}{\end{cases}}
\newcommand{\mn}{\par\medskip\noindent}

\newcommand{\mc}{\mathcal}

\newcommand{\cC}{\mathcal C}

\newcommand{\cK}{\mathcal K}

\newcommand{\cT}{\mathcal T}

\def\Cz{\mathbb{C}}

\def\Nz{\mathbb{N}}

\def\Zz{\mathbb{Z}}

\def\1z{\mathbb{1}}
\newcommand{\lori}{\longrightarrow}
\newcommand{\lole}{\longleftarrow}

\def\SEMI{\mbox{$\times\kern-2pt\vrule height5pt width.6pt \kern3pt $}}

\newcommand{\comment}[1]{}  

\newcommand*{\id}{\mathrm{id}}

\begin{document}

\title[]{The image of Bott periodicity in cyclic homology}
\author{Joachim Cuntz}
\email{cuntz@uni-muenster.de}
\address{Mathematisches Institut, Einsteinstr. 62, 48149 M\"unster, Germany}
\date{\today}

\begin{abstract}
We analyze the relationship between Bott periodicity in topological $K$-theory and the natural periodicity of cyclic homology. This is a basis for understanding the multiplicativity, in odd dimensions, of a bivariant Chern-Connes character from $K$-theory to cyclic theory.
\end{abstract}
\thanks{Gef\"ordert durch Deutsche Forschungsgemeinschaft (DFG) im Rahmen der Exzellenzstrategie des Bundes und der L\"ander EXC 2044 –390685587, Mathematik M\"unster: Dynamik–Geometrie–Struktur}
\maketitle

\section{Introduction}
There are two principal tools for studying the `algebraic topology' of an algebra over $\Cz$: topological $K$-theory on the one hand and periodic cyclic homology (which generalizes de Rham theory to the not necessarily commutative setting) on the other.  One of the first striking observations is the fact that both theories are periodic - they are naturally defined as $\Zz/2$-graded theories. In fact, periodic cyclic theory $HP_*$ has a completely intrinsic $\Zz/2$-periodicity. On the other hand, on topological $K$-theory (for algebras over $\Cz$), the periodicity is due to the non-trivial Bott-periodicity isomorphism.
In the case where the ground field is $\Cz$ and on a natural category of algebras these theories can be compared via a Chern-Connes character (which generalizes the classical Chern character in the commutative case). It is an obvious problem to understand the connection, under the character, between the periodicity operations which are given on the one hand by a `Bott map' and on the other essentially by Connes' $S$-operator.

For such a study, we have to fix a suitable category of topological algebras on which both theories have reasonable properties and can be compared - and where the Bott isomorphism holds. There is such a natural category - the category of locally convex algebras.
Second, for our study we have to choose the correct setting. This is the bivariant one. Both theories have been developed to bivariant theories with an associative composition product. These give a particularly efficient tool for understanding and computing $K$-theoretic and cyclic invariants. Bott periodicity is then naturally defined as a product with a bivariant element. As we will see, a Chern-Connes character can be defined in the bivariant setting. It gives a multiplicative transformation from bivariant $K$-theory, which is denoted on the category of locally convex algebras by $kk_*(A,B)$, to bivariant periodic cyclic theory $HP_*(A,B)$. One immediate problem then consists in determining the image of the bivariant Bott element.

The existence of a bivariant Chern-Connes character in the even case as a multiplicative transformation $kk_0(A,B) \to HP_0(A,B)$ follows directly from the fact that $kk_0$ is the universal functor from the category of locally convex algebras to an additive category subject to three natural properties (half exactness, invariance under smooth homotopies and stability) together with the fact that $HP_*$ shares these properties.
There is an obvious candidate for an extension to the odd case. But unfortunately, for $kk_*$ made $\Zz/2$-periodic, this candidate will not be compatible with the product of two odd elements, because the product uses Bott periodicity. Therefore the extension of the Chern-Connes character from a multiplicative transformation $kk_0 \to HP_0$ to a multiplicative transformation $kk_* \to HP_*$ is linked to a precise understanding of the image of the Bott periodicity isomorphism in $HP_*$.

A key observation in comparing Bott periodicity in $K$-theory to the periodicity in $HP_*$ is the appearance of the numerical factor $2\pi i$. Already in his foundational paper \cite{CoNCG} Connes (motivated by the example of a complex torus and indirectly also by Bott periodicity) used this factor in his definition of the $S$-operator. The factor $2\pi i$ then also appeared in discussions of a Chern-Connes character in \cite{Kar}, \cite{ENN}, \cite{Nist}, \cite{PuAs}, \cite{CuDoc}, \cite{CuWeyl}. As we will see it comes from the connection between the algebra of functions on $[0,1]$, that vanish at 0 and 1, annd the algebra of functions on $S^1$ that vanish in 1, see \cite{Nist}, \cite{CuDoc}. The implication for the behaviour of the product of odd elements in the bivariant theory under the bivariant Chern-Connes character has been discussed in \cite{PuAs}, \cite{CuDoc}, \cite{CuEnz}.

The connection between Bott periodicity in topological $K$-theory and in periodic cyclic homology has been treated in the literature, notably in \cite{Nist}, \cite{CuDoc}. In \cite{Nist} a character has been established for certain bivariant $kk$-elements which is partially multiplicative. In \cite{CuDoc} a fully multiplicative character has been constructed. In both cases the behaviour of Bott periodicity under the character was important. We feel that the issue deserves a more complete presentation. This is the purpose of the present note. An outcome of our discussion will be that a multiplicative character in odd dimensions cannot be entirely canonical, but depends on the choice of a square root of $2\pi i$.

We start by recalling the basics about bivariant $K$-theory $kk_*$. For ease of exposition we use the approach from \cite{CuWeyl}, but applied to the natural class of complete locally convex algebras with submultiplicative seminorms ($m$-algebras). As a first step towards defining the Chern-Connes character we define a multiplicative transformation $\chi_*: kk_*\to HP_*$ from \emph{$\Zz$-graded} (!) bivariant $K$-theory. This transformation is also compatible with the boundary maps in long exact sequences but it does not give a multiplicative transformation from $kk_*$ when we make it $\Zz/2$-periodic using Bott periodicity.

We then give a treatment, which we think is somewhat more systematic and complete than the one in \cite{CuDoc}, of the behaviour of Bott periodicity and of the product under the character and use it to define the final form of the multiplicative Chern-Connes character $ch_*$. This character descends to a transformation from $\Zz/2$-periodic bivariant $K$-theory to $HP_*$. Since the boundary maps in the long exact sequences associated with an extension of algebras are given as products by elements of odd degree, it turns out that $ch_*$ has to multiply the boundary maps by $\sqrt{2\pi i}$. One outcome of our discussion also is that we can alternatively obtain $ch_*$ in the same natural way as $\chi_*$ if we modify the boundary maps in the complex defining $HP_*$ by multiplying them by $\sqrt{2\pi i}$. Denote the homology of the so modifies complex by $\overline{HP}_*$. Then $ch_*$ will be compatible with the boundary maps in $\overline{HP}_*$, see Remark \ref{bound}(b).

We mention though that the approach we describe here is fundamentally equivalent to the one in \cite{CuDoc}. We also note that our arguments carry over without a problem to arbitrary locally convex algebras if we use the framework of \cite{CuTh}.

\section{The setting}\label{set}
For convenience we will work with the simplest natural category suitable for the discussion of the image of the Bott map. This is the category of algebras over $\Cz$ with a complete locally convex topology defined by a family of submultiplicative seminorms. Following \cite{CuDoc} we call these algebras $m$-algebras. The natural tensor product in this category is the completed projective tensor product, denoted by $A\otimes B$. By an extension of $m$-algebras we mean an exact sequence $0\to I\stackrel{j}{\lori} E\stackrel{q}{\lori} B\to 0$ of $m$-algebras where $j,q$ are continuous homomorphisms. Moreover we will always assume that an extension admits a continuous linear split, i.e. a continuous linear map $s: B\to E$ such that $q\circ s=\id$. Using shorthand we write $I\to E \to B$ for such an extension.

\emph{The universal extension.} Let $A$ be an $m$-algebra. The non-unital tensor algebra
$T_{alg}A = A\oplus A^{\otimes 2}\oplus A^{\otimes 3}\oplus \cdots$ can
be completed to an $m$-algebra $TA$ such that this becomes the universal $m$-algebra generated by an image of $A$ under a continuous linear map (the linear inclusion of $A$ onto the first direct summand in $TA$).
\\There is a natural
surjective homomorphism $TA\to A$ mapping a tensor of the form $a_1\otimes \ldots \otimes a_n$ to the product $a_1 \ldots a_n$ in $A$. We denote its kernel by
$JA$. The extension $JA \to TA \to A$ is universal among all extension of $A$ admitting a continuous linear splitting.
Indeed, if $I \to E\to A$ is any extension of $m$-algebras
with a continuous linear splitting $s: A\to E$, then there is a
canonical commutative diagram of the form
$$
\begin{array}{ccccccccc} 0 &\to &  J A & \to & T A  & \to &
A
& \to & 0\\[-0cm]
&      &  \quad\downarrow {\gamma_s}& & \downarrow & &
\quad\downarrow\id
 & &\\[-0cm] 0
& \to &  I & \to &  E & \to &  A & \to & 0
\end{array}
$$
where the map $TA \to E$ is the unique algebra homomorphism extending the linear map $TA \supset A \stackrel{s}{\lori} E$ and $\gamma_s$ its restriction to $JA$. We call the map $\gamma_s : JA \to I$ the classifying map for the extension. It depends on $s$ only up to smooth homotopy, see \cite{CuDoc} for details.

We also have to specify the class of functions that we allow for homotopies. Since $HP_*$ is only invariant under differentiable homotopies, this will be the class of infinitely differentiable, or smooth, functions on intervals. Given an $m$-algebra $A$, we denote by $A[0,1], A(0,1], A(0,1)$ the algebras of smooth $A$-valued functions on $[0,1]$ with all higher derivatives vanishing at $0$ and $1$, the functions themselves vanishing at $0$ in the second case and at $0$ and $1$ in the third case.
Two homomorphisms $\varphi_0,\varphi_1: A\to B$ are called diffotopic (or differentiably homotopic), if there is a homomorphism $\Phi :A \to B[0,1]$ such that $ev_0\circ\Phi=\varphi_0$ $ev_1\circ\Phi = \varphi_1$.\\
Finally, we need an appropriate notion of stabilization of an $m$-algebra $A$ by matrices. The simplest choice for our purposes here is given by the completed projective tensor product $\cK\otimes A$ of $A$ by the algebra $\cK$ of $\Nz\times\Nz$-matrices with rapidly decreasing coefficients (which is a natural $m$-algebra completion of the algebra $M_\infty$ of finite matrices of arbitrary size, see \cite{CuDoc}).

We can now recall the definition, from \cite{CuWeyl}, of the bivariant theory $kk_*$ slightly adapted to the setting we are considering here.
Denote by $\gamma_A : JA \to A(0,1)$ the classifying map for the cone extension $A(0,1)\to A(0,1]\to A$. We denote by $[A,B]$ the set of diffotopy classes of homomorphisms $A\to B$. For an $m$-algebra $D$, by we define $J^nD$ recursively by $J^{n+1}D = J(J^nD)$. Using the map $\gamma_A: JA \to A(0,1)$ for each $m,n$ we get a map $[J^nA,B(0,1)^m]\lori [J^{n+1}A,B(0,1)^{m+1}]$.\\
Given two $m$-algebras $A,B$ we can now define
$$kk_n(A,B)=  \lim_{\mathop{\lori}\limits_{k}} [J^{k-n}A,\cK\otimes B(0,1)^k ] $$
(note that $n<0$ is allowed).
\\$kk_n$ is an abelian group with the addition
$$[X,\cK\otimes Y]\times [X,\cK\otimes Y] \to [X,\cK\otimes Y]$$
induced for $m$-algebras $X,Y$ by the inclusion $\cK\oplus \cK\hookrightarrow \cK$

There is a
natural product $kk_s(A_1,A_2)\times kk_t(A_2,A_3)\lori kk_{s+t}(A_1,A_3)$ For simplicity of notation we explain it for $s=t=0$. The general case follows from the identities $kk_s(A,B)=kk_0(A,B(0,1)^s)$, $kk_{-s}(A,B)=kk_0(J^s A,B)$.\\
Assume that two elements $a$ in $kk_0(A_1,A_2)$ and $b$ in $kk_0(A_2,A_3)$ are
represented by $\alpha :J^nA_1\to \cK\otimes A_2(0,1)^{n}$ and
$\beta :J^mA_2\to \cK\otimes A_3(0,1)^{m}$, respectively. We
define their product $a\cdot b$ in $kk_{0}(A_1,A_3)$ by
the following composition of maps
\bglnoz J^m(J^nA_1)\quad \mathop{\lori}\limits^{J^m(\alpha)}\quad
J^m(\cK\otimes A_2(0,1)^{n}) \,\to\qquad\qquad\qquad\qquad \\ \qquad\qquad\cK\otimes (J^m A_2)(0,1)^{n}
\quad \mathop{\lori} \limits^{\beta(0,1)^{n}}\quad \cK\otimes
A_3(0,1)^{n+m}\eglnoz
The second arrow is defined as the
composition of the natural map $J^m(A_2\otimes E)\to J^m(A_2)\otimes E$
(choosing $E=\cK\otimes (0,1)^n$), that exists, by definition of $J$, for any $E$, with the isomorphism of $\cK\otimes (J^m A_2)(0,1)^n$ that switches the orientation of one
of the intervals $(0,1)$ in $(0,1)^n$ according to the parity
of $mn$ (i.e. we take $(-1)^{mn}$ times the diffotopy class
of this map). Of course we also use the isomorphism $\cK\otimes\cK\cong\cK$. The sign is necessary because the map $$J^m(\cK\otimes A_2(0,1)^n) \,\to\, \cK\otimes (J^m A_2)(0,1)^n$$
permutes the noncommutative $m$-fold suspension $J^m$ with the
$n$-fold ordinary suspension $A(0,1)^n$. The sign is explained by the following easy lemma that then shows that the product is
well defined (i.e. does not depend on the choice of
representatives for $a$ and $b$ in the inductive limits
defining $kk_0$).

\blemma\label{well} For
an $m$-algebra $D$ as above we denote by $\gamma_D:JD\to
D(0,1)$ the classifying map for the cone extension
$D(0,1)\to D(0,1]\to D$.\\ (a) Let
$I\to E \mathop{\lori} \limits^{\pi} A$ be a
linearly split extension of locally convex algebras and
$\alpha$ its classifying map. Then the two maps $J^2A\to
I(0,1)$ defined by the compositions
$J^2A\mathop{\lori}\limits^{J(\alpha)}JI
\mathop{\lori}\limits^{\gamma_I} I(0,1)$ and
$J^2A\mathop{\lori}\limits^{-\gamma_{JA}} (JA)(0,1)
\mathop{\lori}\limits^{\alpha(0,1)}I(0,1)$ are
diffotopic. Here $-\gamma_{JA}$
stands for $\gamma_{JA}$ followed by a change of orientation on $(0,1)$. \\
(b) The natural map $\gamma_{JA}:\,J^2A\to (JA)(0,1)$ is
diffotopic to the composition of the natural maps $J^2A\to J(A(0,1))$ and $J(A(0,1))\to
(JA)(0,1)$. \elemma
\bproof (a) The two compositions are
classifying maps for the two 2-step extensions in the first
and last row of the following commutative diagram
\[
\begin{array}{ccccccccccc} 0 &\to &  I(0,1) & \to & I(0,1]  &
\to &   E&\to &A
& \to & 0\\[0.1cm]
&      &  \parallel& &
\downarrow & & \downarrow & &\downarrow& &\\[0.1cm] 0
& \to &  I(0,1) & \to &  E(0,1] & \to &  Z_\pi & \to & A&\to
& 0\\[0.1cm]&      &  \parallel& &
\uparrow & & \uparrow & &\uparrow& &\\[0.1cm]
0 & \to &  I(0,1) & \to &  E(0,1) & \to &  A[0,1) & \to &
A&\to & 0
\end{array}
\]
where $Z_\pi\subset E\oplus A[0,1]$ is the mapping cylinder.
They are diffotopic since they also are classifying maps for
the extension in the middle row. Here we use the
cone $A[0,1)$ in place of $A(0,1]$ and thus reverse
the orientation of $[0,1]$.\\
(b) The first map is classifying for the first row in
the commutative diagram \[
\begin{array}{ccccccccccc} 0 &\to &  (JA)(0,1) & \to & (TA)(0,1)
& \to &   A(0,1]&\to &A
& \to & 0\\[0.1cm]
&      &  \quad\uparrow{\scriptstyle\varphi} & &
\uparrow & & \uparrow& & \uparrow & &\\[0.1cm] 0
&\to &  J(A(0,1)) & \to & T(A(0,1)) & \to &   A(0,1]&\to &A &
\to & 0
\end{array}\]
while the second map is the composition of the classifying
map for the second row composed with the first vertical arrow
$\varphi$ in the diagram - thus also classifying for the
first row. \eproof
Note that every continuous homomorphism
$\alpha:A\to B$ induces an element $kk(\alpha)$ in $kk_0(A,B)$ and that $kk(\beta\circ\alpha)=kk(\alpha)\cdot
kk(\beta)$. We denote by $1_A$ the element $kk(\id_A)$
induced by the identity homomorphism of $A$. It is a unit in
the ring $kk_0(A,A)$.\\
From its definition one derives the following properties of the functor $kk$.
\begin{itemize}
  \item[(E1)]  Differentiable homotopy invariance: The evaluation maps $A[0,1]\to A$ induce $kk_0$-equivalences.
  \item[(E2)]  Stability: The natural inclusion $A\hookrightarrow \cK\otimes A$ induces a $kk_0$-equivalence.
  \item[(E3)]  Half-exactness: Every extension $I\to E\to A$ with a continuous linear split induces short exact sequences in both variables of $kk_0$.
\end{itemize}
\bremark By property (E3) combined with (E1) every extension $I\to E\to A$ of $m$-algebras induces long exact sequences in both variables of $kk_*$. The boundary maps $kk_n\to kk_{n-1}$ are given by product with the element $kk(\gamma) \in kk_0(JE,I)=kk_{-1}(E,I)$ for the classifying map $\gamma$ of the extension.\eremark
We know that $kk_0$ is best viewed as an additive category with objects $m$-algebras and morphisms $kk_0(A,B)$. Then $\alpha\to kk(\alpha )$ defines a functor from $m$-algebras to the category $kk_0$ which is universal with the properties above:
\btheo\label{uni}
Let $E$ be a functor from $m$-algebras to an additive category $\cC$ which satisfies the properties $E1,E2,E3$ in both variables. Then there is a unique functor $F: kk_0 \to \cC$ such that $E = F\circ kk$.\etheo
\bproof
We sketch the simple essence of the argument here and give a detailed construction of $F$ (for the case $E=HP$) in section \ref{BiCC}.\\ The properties of half-exactness and homotopy invariance imply that any extension $I\to A \to B$ of $m$-algebras with a continuous linear split induces long exact sequences under $E$ in both variables (in particular it follows that $E$ is automatically additive). We apply the long exact sequences to the standard extensions
$$JA\to TA\to A \qquad
A(0,1)\to A(0,1]\to A$$
Write $E(A,B)$ for the image, under $E$ of the $m$-algebra morphisms $A\to B$ in the morphisms of $\cC$. The fact that $E$ applied to the contractible algebras $TA$ and $A(0,1]$ gives 0, shows that there are isomorphisms
$E(JA,B(0,1))\cong E(A,B)$ and by iteration $E(J^nA,B(0,1)^n)\cong E(A,B)$. Using stability of $E$ we get moreover $E(J^nA,\cK \otimes B(0,1)^n)\cong E(A,B)$. Let $x$ be an element of $kk_0(A,B)$ represented by a homomorphism $\varphi :J^n A  \to \cK \otimes B(0,1)^n$. We set $F(x)= E(\varphi)\in E(J^nA,\cK\otimes B(0,1)^n)\cong E(A,B)$. The explicit description of $F$ in section \ref{BiCC} will make it obvious that $F$ is well defined and a functor (i.e. multiplicative).
\eproof
We also note that $kk_0(\Cz,B)$ gives the usual topological $K$-theory (e.g. when applied to a Banach algebra or to a Fr\'{e}chet algebra), see \cite{CuDoc}, Section 7.

\section{Bott periodicity}\label{BP}
The well known proofs of Bott periodicity, e.g. \cite{AtB}, \cite{KasInv}, \cite{CuKth} use explicitly or implicitly an index map. There is a universal model for the index map that is encoded in the Toeplitz extension. Algebraically, the Toeplitz algebra (for our purposes over the field $\Cz$) is the universal unital $\Cz$-algebra $\cT$ generated by two elements $u,v$ satisfying the relation $vu=1$. It is easy to check that the ideal generated by the defect projection $e=1-uv$ in $\cT$ is isomorphic to the algebra $M_\infty$ of finite matrices of arbitrary size, and the quotient by this ideal is isomorphic to the algebra $\Cz[z,z^{-1}]$ of Laurent polynomials. There is a well known simple algebraic argument that goes back to \cite{CuKth} showing that the natural inclusion $\Cz \to \cT$ induces an isomorphism $E(\Cz) \to E(\cT)$ for every homotopy invariant (under polynomial homotopies), split-exact and matrix stable functor $E$. As a consequence, if we denote by $\cT_0$ the kernel of the natural map $\cT\to \Cz$, then $E(\cT_0)=0$. The argument for these facts is so simple and general that it carries over to all reasonable completions of $\cT$, in particular to natural C*-algebra completions or $m$-algebra completions. The purely algebraic form of the argument has been used more recently also in \cite{CoTh}.

For our purposes we use $m$-algebra completions $\overline{\cT}, \overline{\cT}_0$ of $\cT,\cT_0$, see \cite{CuDoc}, \cite{CuEnz}. We get extensions
$\cK\to \overline{\cT} \to C^\infty(S^1)$,  $\cK\to \overline{\cT}_0 \to C^\infty(S^1\setminus 1)$ where $C^\infty(S^1\setminus 1)$ denotes the algebra of $C^\infty$-functions on the circle $S^1$ that vanish at the point 1. The general argument mentioned above then shows that $\overline{\cT}_0$ is equivalent to 0 in $kk_0$ see \cite{CuDoc}. This means in particular that $kk_n(\overline{\cT}_0\otimes A,B) $ and $kk_n(A,\overline{\cT}_0\otimes B)$ are 0 for all $m$-algebras $A,B$ and all $n$.

Bott periodicity for $kk_*$ means that $$kk_{n+2}(B,A) =kk_n(B,A(0,1)^2) \cong kk_n(B,A)$$ or, equivalently, that $kk_{n-2}(A,B) =kk_n(J^2A,B)\cong kk_n(A,B)$ for all $A,B,n$. These periodicity isomorphisms are implemented by an invertible element in $kk_{-2}(A,A)$ or by its inverse in $kk_2(A,A)$.
This element has a natural representation. We note that any extension $I\to D \to Q$ of $m$-algebras with a continuous linear split where $D$ is contractible, or more generally equivalent to 0 in $kk_0$, defines an invertible element in $kk_{-1}(Q,I)$ (this follows from the long exact sequences in both variables). Now, given an $m$-algebra $A$, we have two extensions of that type
$$\cK\otimes A\to \overline{\cT}_0\otimes A \to C^\infty(S^1\setminus 1)\otimes A\qquad A(0,1)\to A(0,1]\to A$$
They determine invertible elements $\tau_A$ in $kk_{-1}(C^\infty (S^1\setminus 1)\otimes A,\cK\otimes A)$ and $\rho_A$ in $kk_{-1}(A,A(0,1))$, respectively. Moreover the natural inclusion $\Cz(0,1) \to C^\infty (S^1\setminus 1)$ as functions on $S^1$ vanishing with all derivatives at 1, determines an element $\sigma_A$ in $kk_0(A(0,1), C^\infty (S^1\setminus 1) \otimes A)$. It is invertible by differentiable homotopy invariance of $kk_0$ (by an argument as e.g. in \cite{CuDoc}, 1.1, p.143 there is a homotopy inverse to $\sigma$). Finally denote by $\kappa_A$ the element in $kk_0(A,\cK\otimes A)$ induced by the natural inclusion $A \to \cK\otimes A$.
The product $\rho_A\sigma_A\tau_A\kappa_A^{-1}$ defines an invertible element in $kk_{-2}(A,A)$. This fundamental element has been used explicitly or implicitly in many places, e.g. \cite{CuKth}, \cite{KasInv}, \cite{Nist}, \cite{CuDoc}. We record it in the following definition
\bdefin\label{Bott} For an $m$-algebra $A$ we denote by $b_A$ the element $\rho_A\sigma_A\tau_A\kappa_A^{-1}$ of $kk_{-2}(A,A)$ as above. \edefin

\section{A multiplicative transformation $kk_*\to HP_*$}\label{BiCC}
We now come to the connection of $kk_*$ to bivariant periodic cyclic homology $HP_*$. Of course we have to use the version of $HP_*$ for topological algebras. This means that for $m$-algebras, in the algebraic definition \cite{CQCNs} of $HP_*$ we have to replace algebraic tensor products by completed projective tensor products.
Recall then that any extension $I\to D\to A$ of $m$-algebras with a continuous linear splitting induces long exact sequences in both variables of $kk_*$. For such an extension the classifying map $\gamma :JA \to I$ gives an element in $kk_{-1}(A,I)$. It is important for our purposes to note that the boundary maps $kk_n\to kk_{n-1}$ in both long exact sequences under $kk_*$ are given by the product by this element.

By excision \cite{CQEx},\cite{CuExTop} $HP$ satisfies (E3). It therefore has long exact sequences in both variables. Given an $m$-algebra $A$ and another $m$-algebra $B$ we consider the long exact sequences in $HP_*(B,A)$ induced in the second variable by the universal extension $JA\to TA\to A$ on the one hand and by the cone extension $A(0,1)\to A(0,1]\to A$ on the other hand. We denote by $\delta$ the boundary map for the first sequence, by $\bar{\delta}$ the boundary maps for the second sequence and obtain a commutative diagram \\[-0.3mm]
\bgl\label{boundary}
\begin{array}{ccccccc}  &\to &  HP_{0}(B,A) & \stackrel{\delta}{\lori} & HP_1(B,JA)  & \to &
\\[-0cm]
&      &  \quad\downarrow\, = & & \downarrow {\cdot\,\gamma_A} & &
\\[-0cm]
& \to &  HP_0(B,A) & \stackrel{\bar{\delta}}{\lori} &  HP_1(B,A(0,1)) & \to &
\end{array}
\egl

Here $\gamma_A: JA \to A(0,1)$ denotes the natural map and $\cdot \,\gamma_A$ denotes product by $\gamma_A$ on the right (note that here and elsewhere we denote the product in $kk_*$ and in $HP_*$ in the order opposite to the composition of homomorphisms). Since $HP_*(B,TA)=HP_*(B,A(0,1])=0$ by contractibility of $TA$ and $A(0,1]$, we see that $\delta, \bar{\delta}$ are isomorphisms.

Thus there are invertible elements $\alpha_A\in HP_1(A,JA)$ and $\beta_A\in HP_1(A(0,1),A)$ such that $\delta=\cdot\,\alpha_A$ and $\bar{\delta}=\cdot\,\beta_A^{-1}$ (explicitly they are determined by $\alpha_A=\delta(1_A)$, $\beta_A^{-1}=\bar{\delta}(1_A)$, see \cite{CQEx} Theorem 5.5). The commutative diagram \eqref{boundary} shows that $\alpha_A\gamma_A\beta_A = 1_A$ where
$1_A$ denotes the identity element of $HP_0(A,A)$ and also $\gamma_A\cdot\beta_A\cdot\alpha_A =1_{JA}$, $\beta_A\cdot\alpha_A\cdot\gamma_A =1_{A(0,1)}$.\\
By iteration we get
invertible elements $\alpha_A^n= \alpha_A\alpha_{JA}\,\cdots\, \alpha_{J^nA}\in$ $HP_n(A,J^nA)$ and
$\beta_A^n=\beta_{A(0,1)^n}\,\cdots\, \beta_{A(0,1)}\beta_A\in HP_n(A(0,1)^n,A)$. Let moreover $\gamma^n$ be
the element of $HP_0(J^nA,A(0,1)^n)$ induced by the canonical
homomorphism $J^nA\to A(0,1)^n$. The following relations are derived by iteration from diagram \eqref{boundary}
\bgl\label{gamma}\alpha_A^n\cdot\gamma^n\cdot\beta_A^n =1_A\qquad
\gamma^n\cdot\beta_A^n\cdot\alpha_A^n =1_{J^nA}\qquad
\beta_A^n\cdot\alpha_A^n\cdot\gamma^n =1_{A(0,1)}\egl
As a first step towards constructing a multiplicative Chern-Connes character we can now define a multiplicative transformation $\chi$ from the $\Zz$-graded theory $kk_*$ to $HP_*$.
\btheo\label{Chern} There are
natural maps $\chi_s:kk_s(A,B)\to HP_s(A,B)$ such that
$$\chi_{s+t}(\varphi\cdot\psi)=\chi_s(\varphi)\cdot \chi_t(\psi)$$ in
$HP_{s+t}(A,C)$ for $\varphi\in kk_s(A,B)$ and
$\psi\in kk_t(B,C)$ and such that $\chi_0(1_A)=1_A$ for
each $A$. Here $s$ and $t$ are in $\Zz$ and $HP_s$ is
2-periodic in $s$.\etheo
\bproof Ignoring
first the stabilization by $\mc K$ assume that $\varphi\in
kk_s(A,B)$ is represented by a homomorphism $\eta
:J^{n-s}A\to B(0,1)^n$. We then define
$$\chi_s(\varphi)= \alpha_A^{n-s}HP_0(\eta)\beta_B^n$$
The relation $\alpha_A^n\cdot\gamma^n\cdot\beta_A^n =1_A$
generalizes to
$\alpha_A^n\cdot\theta^n\cdot\beta_B^n =\theta$ if $\theta$
is an element of $HP_0(A,B)$ determined by a homomorphism
$A\to B$ and $\theta^n$ is the element of
$HP_0(J^nA,B(0,1)^n)$ determined by the homomorphism $J^nA\to
B(0,1)^n$ induced by $\theta$. Since the inductive limit in
the definition of $kk_s$ is taken exactly with
respect to the passage from homomorphisms $A\to B$ to
homomorphisms $J^nA\to B(0,1)^n$, this shows that $\chi$ is
well defined. The compatibility with the product follows from the relations (\ref{gamma}) and
Lemma \ref{well}.\\
Finally, the natural inclusion $B\to \mc
K\otimes B$ induces an invertible element $\iota$ in $HP_0(B,\mc
K\otimes B)$. The formula for $\chi$ then becomes
\bgl\label{forchi}\chi_s(\varphi)=
\alpha_A^{n-s}HP_0(\eta)\beta_B^n\,\iota^{-1}\egl
if $\varphi$
is represented by the map $\eta :J^{n-s}A\to\mc K\otimes
B(0,1)^n$. The stabilization by $\mc K$ does not lead to any
essential change in the arguments for well-definedness and
multiplicativity. \eproof
\bremark Note that the formula for $\chi_0$ makes the universal functor $kk_0\to HP_0$ from Theorem \ref{uni} explicit.\eremark
\section{The image of Bott periodicity in $HP_*$ and the final form of the Chern-Connes character}
In section \ref{BiCC} we have seen that, essentially due to the universal property of the $kk$-functor, there is a completely canonical multiplicative transformation $\chi$ from the $\Zz$-graded(!) theory $kk_*$ to the a priori also $\Zz$-graded theory $HP_*$. But $HP_*$ is actually $\Zz/2$-periodic in a canonical way (due to Connes' $S$-operator). On the other hand $kk_*$ is periodic too - not by definition but using the product with the Bott element $b_A$ from Definition \ref{Bott}.

The following computation has been made first in \cite{Nist}.
\btheo\label{imbott}
(a) One has $\chi_{-2}(b_\Cz) = (2\pi i)^{-1}1_\Cz$.\\
(b) Given an $m$-algebra $A$ one has  $\chi_{-2}(b_A) = (2\pi i)^{-1}1_A$.
\etheo
\bproof
(a) The element $b_\Cz$ is defined as the product $\rho_\Cz\sigma_\Cz\tau_\Cz\kappa_\Cz^{-1}$ where $\tau_\Cz$, $\rho_\Cz$ are the elements of $kk_{-1}$ determined by the extensions
\bgl\label{ext}\cK\to \overline{\cT}_0\to C^\infty(S^1\setminus 1)\qquad \Cz(0,1)\to \Cz(0,1]\to \Cz\egl
and $\sigma_\Cz,\kappa_\Cz$ are the isomorphisms in $kk_0$ induced by the maps $\Cz(0,1) \to C^\infty(S^1\setminus 1)$ and $\Cz\to \cK$. By naturality of the boundary maps in $kk_*$, $b_\Cz$ can also be written as the product $\rho_\Cz\overline{\sigma}_\Cz\overline{\tau}_\Cz\kappa_\Cz^{-1}$ where $\overline{\tau}_\Cz$ denotes the element of $kk_{-1}$ determined by the extension $\cK\to \overline{\cT}\to C^\infty(S^1)$ and $\overline{\sigma}_\Cz$ the element of $kk_0$ induced by the inclusion $\Cz(0,1)\to C^\infty(S^1)$.

Given any extension $I\to E\to Q$ with classifying map $\gamma:JQ \to I$, the product by $\chi_{-1} (\gamma)$ gives the boundary maps in the long exact sequences associated by the extension in the first or second variable of $HP$, respectively, see \cite{CQEx}, Theorem 5.5. So the product by $\chi_{-1}(\rho_\Cz)$ and $\chi_{-1}(\overline{\tau}_\Cz)$ on the right gives the boundary maps $HP_*(\Cz,\Cz) \to HP_{*-1}(\Cz,\Cz(0,1))$ and $HP_*(\Cz,C^\infty(S^1)\to HP_{*-1}(\Cz, \cK)$ for the extensions $\Cz(0,1)\to \Cz(0,1]\to \Cz$ and $\cK\to \overline{\cT}\to C^\infty(S^1)$, respectively. We are going to determine these boundary maps using that $HP_*(\Cz,A)=HP_*(A)$.\\
In fact, since all the algebras involved have small homological dimension, it turns out that we can do the computation in $HX_*$ rather than in $HP_*$. Recall that $HX_*(A)$ denotes the homology of the $\Zz/2$-graded $X$-complex
$$
X(A):\qquad A\qquad
\begin{array}{c}
\mathop{\lori}\limits^{\natural d} \\[-6pt]
\mathop{\lole}\limits_b
\end{array}
\qquad \Omega^1(A)_{\natural}
$$
which is the quotient of the periodic $B-b$-complex $\widehat{\Omega}_* (A)$ by the first subcomplex $F_1$ in the Hodge filtration. Here $\natural$ denotes quotient by commutators. One has \bgl\label{HX}HX_0(A)=SHC_2(A) \qquad HX_1(A) = HC_1(A)\egl where $HC_*$ denotes ordinary cyclic homology, see \cite{CQCNs}, p.393.\\
By the argument in \cite{CoNCG}, p.127 and formula (4) above one has $HC_1(\Cz(0,1))=HX_1(\Cz(0,1))=\Cz$, $SHC_0(\Cz(0,1))=HX_0(\Cz(0,1))=0$.  Using the fact that $\Cz(0,1)$ is $h$-unital and excision for $HC_*$ one gets that $HX_*(\Cz(0,1])=0$. For the extension $\Cz(0,1)\to \Cz(0,1]\to \Cz$ we then get the following commutative diagram
\bgl\label{cone}
 \begin{array}{ccccccc}
 \lori & HP_*(\Cz(0,1)) & \lori & HP_*(\Cz(0,1]) & \lori &
 HP_*(\Cz) & \lori  \\[2pt]
 & \:\downarrow \cong && \downarrow \cong &
 & \;\downarrow \cong &
\\[3pt]
 \lori & HX_*(\Cz(0,1):\Cz(0,1]) & \lori &
 HX_*(\Cz(0,1]) & \lori &
 HX_*(\Cz) & \lori
 \end{array}
\egl
 where $HX_*(\Cz(0,1):\Cz(0,1]$ denotes the homology of the relative complex i.e. the homology of the kernel of $X_*(\Cz(0,1])\to X_*(\Cz)$.
The rows are exact and the first downward arrow is an isomorphism by the five-lemma. Since $\Cz(0,1)$ is $h$-unital $HX_*(\Cz(0,1):\Cz(0,1])$ is isomorphic to $HX_*(\Cz(0,1))$.\\
Similarly, we get a commutative diagram for the Toeplitz extension\mn
\bgl\label{Toep}
 \begin{array}{ccccccc}
 \lori & HP_*(\cK) & \lori & HP_*(\overline{\cT}) & \lori &
 HP_*(C^\infty(S^1)) & \lori  \\[2pt]
 & \:\downarrow \cong && \downarrow \cong &
 & \;\downarrow \cong &
\\[3pt]
 \lori & HX_*(\cK) & \lori &
 HX_*(\overline{\cT}) & \lori &
 HX_*(C^\infty(S^1)) & \lori
 \end{array}
\egl
The second row is exact because $\cK$ is $h$-unital. The first and third downward arrows are isomorphisms because $\cK$ and $C^\infty(S^1)$ have homological dimension 1, and finally the middle downward arrow is an isomorphism by the five-lemma.\\
We now have to identify the boundary maps $\delta_1, \delta_2$ for the second rows in \eqref{cone} and \eqref{Toep}. They will represent $\chi_{-1}(\rho_\Cz)$ and  $\chi_{-1}(\overline{\tau}_\Cz)$. By definition of the boundary map in the long exact sequence for an extension of complexes, in \eqref{cone} we have to lift the element 1 in $X_0(\Cz)$ to an element in $X_0(\Cz(0,1])$. For such a lift we choose a monotone function $h\in \Cz(0,1]$ such that $h(1)=1$. We then apply the boundary operator $\natural d$ in the $X$-complex. This gives $\delta_1(1)=\natural dh$. \\
Under the map $X(\Cz(0,1))\to X(C^\infty (S^1))$ induced by the inclusion $\overline{\sigma}:\Cz(0,1)\to C^\infty (S^1)$, $\natural dh$ is mapped to $(2\pi i)^{-1} w^{-1}dw$ where $w=e^{2\pi ih}$ (this is because modulo commutators $w^{-1}dw=e^{-2\pi ih}\,2\pi i\, dh\,e^{2\pi ih}$). By homotopy invariance of $HP_*$ the element $w^{-1}dw$ becomes in $HP_*(C^\infty (S^1))=HX_*(C^\infty (S^1))$  equal to $z^{-1}dz$ with $z=e^{2\pi it}$ the standard generator of $C^\infty (S^1)$. Now, to compute $\delta_2(z^{-1}dz)$, we lift $z,z^{-1}$ to $u,v$ and $z^{-1}dz$ to $vdu$ in $X_1(\overline{\cT})$. When we apply the boundary $b$ in the $X$-complex to $vdu$, we obtain $vu - uv=e$.
In conclusion $\delta_2 \overline{\sigma}_*\delta_1(1)= (2\pi i)^{-1} e$ with $e$ the standard minimal idempotent in $\cK$ (here $\overline{\sigma}_*$ denotes the map induced on $HX_0$ by $\overline{\sigma}$).\\
(b) follows from (a) by exterior product, see e.g. \cite{CQEx}, p.86. In fact $b_A$ is the exterior product $b_\Cz\otimes \id_A$.\\
\eproof
Since the ideals $\cK$ and $\Cz(0,1)$ are h-unital, the computation above can in fact be done in non-periodic cyclic homology and Nistor has found in \cite{Nist} that $b_A$ corresponds in the non-periodic case to $2\pi i\, S$.
\bremark
In \cite{CuDoc} we had given a proof for the formulas in Theorem \ref{imbott} that reduces the computation of the boundary maps to the purely algebraic case of finitely generated algebras over $\Cz$. For the present paper we have chosen the  proof given above in order to emphasize the role of the inclusion map $\overline{\sigma}$ from $\Cz(0,1)$ into $C^\infty(S^1)$.\eremark
Recall the invertible elements $\alpha_A^n= \alpha_A\alpha_{JA}\,\cdots\, \alpha_{J^nA}\in$ $HP_n(A,J^nA)$ and
$\beta_A^n=\beta_{A(0,1)^n}\,\cdots\, \beta_{A(0,1)}\beta_A\in HP_n(A(0,1)^n,A)$ from section \ref{BiCC}. Left and right multiplication by these elements gives the transformation $\chi_n : kk_n\to HP_n$.
In particular $\chi_0: kk_0(A,B) \to HP_0(A,B)$ is the universal functor described in section \ref{set}. By the discussion in section \ref{BiCC} it sends a homomorphism $\varphi: J^nA \to B(0,1)^n$ representing an element of $kk_0(A,B)$ to $\alpha_A^n \varphi \beta_B^n $. This part of $\chi$ therefore is the canonical choice for the Chern-Connes character $kk_0\to HP_0$.\\ Now what happens to $\chi_{2n}:kk_{2n}\to HP_{2n}$. The group $kk_n (A,B)$ is isomorphic to $kk_{n-2} (A,B)$ under left multiplication $x\mapsto b_A\,\cdot x$ by the Bott element and the two groups are identified under this isomorphism in the periodic $\Zz/2$-graded $K$-theory. We have to check how this identification behaves under $\chi$. This means that we have to compute $\chi_{n-2} (b_A\,\cdot x)$ in terms of $\chi_n(x)$.\\
But this has been achieved in Theorem \ref{imbott} and the answer is $\chi_{n-2} (b_A\,\cdot x)=$ $(2\pi i)^{-1} \,\chi_n(x)$. Therefore, if we define
\bgl\label{ev} ch_{-2n}(x) = (2\pi i)^n \chi_{-2n}(x)\egl
for an element $x\in kk_{-2n}(A,B)$ then we get a multiplicative transformation $ch$ from the even part $kk_{2*}(A,B)$ of $kk$ to $HP_{0}(A,B)$. Moreover
$$ch_{-2n-2}(b_A\cdot x)=(2\pi i)^{n+1}\chi_{-2n-2}(b_A\cdot x)=(2\pi i)^{n+1}(2\pi i)^{-1}\chi_{-2n}(x)=ch_{-2n}(x).$$
Thus $ch_*$ is invariant under $b_A\,\cdot$ and descends to a transformation $kk_{ev}\to HP_{0}$ from $kk_{2\star}$ made $\Zz/2$-periodic. This is the natural choice for the Chern-Connes character in the even case.

The situation is more complicated for the odd part of the character. We have to extend $ch$ from the even part to a transformation $ch_*:kk_*\to HP_*$ that is multiplicative and descends to $kk_*$ made periodic via the Bott map. Multiplicativity means in particular that $ch_{-1}(x)ch_{-1}(y)=ch_{-2}(xy) =2\pi\,i\chi_{-2}(xy)$ for $x,y\in kk_{-1}(A,A)$.
This equation leads more or less inevitably to the following definition.
\bdefin\label{od}
We define the bivariant Chern-Connes character by $ch_n(x) = (2\pi i)^{-n/2}\chi_n(x)$ given $x\in kk_n(A,B)$.
\edefin
Thus, more explicitly, if $x$ in $kk_n(A,B)$ is represented by $\varphi: J^{k-n}A \to B(0,1)^k$, then we set
\bgl\label{chi2} ch_n(x)= (2\pi\,i)^{-n/2}\alpha_A^{k-n}HP_0(\varphi)\beta_B^k\egl
For more clarity, in this formula we have suppressed the identification of $B$ with $\cK\otimes B$ via $\iota^{-1}$ as in \eqref{forchi}. Observe that we cannot avoid the non-canonical choice of a square root for $2\pi i$ that we have to make once and for all.
\btheo
The transformation $ch_*$ defined in \eqref{od} extends the definition \eqref{ev} to the case of odd $n$. It is multiplicative, i.e. it satisfies $ch_s(x)ch_t(y) = ch_{s+t}(xy)$ and it descends to a multiplicative transformation $kk_{ev/od}\to HP_{ev/od}$ from $kk_*$ made periodic via the Bott map.
\etheo
\bproof Multiplicativity follows from the multiplicativity of $\chi$. The compatibility with the Bott map follows as before from Theorem \ref{imbott}.
\eproof
\bremark\label{bound}
(a) Given an extension $I\to A\to B$ of $m$-algebras denote by  $x$ in $kk_{-1}(B,I)$ the element defined by the classifying map $JB\to I$. The boundary maps in the long exact sequences for $kk_*$ and $HP_*$ are given as the product by $x$ and $\chi(x)$ respectively. Since $\chi$ is multiplicative it therefore respects the boundary maps (see also \cite{Nist2} for a more precise statement). On the other hand $ch_{-1}(x)=\sqrt{2\pi i}\chi_{-1}(x)$. Therefore $ch_*$ multiplies the boundary maps by $\sqrt{2\pi i}$.\\
(b) Consider the theory $\overline{HP}_*$ which is obtained if we multiply the boundary map $B-b$ in the complex defining $HP_*$ by $\sqrt{2\pi i}$. Since the kernels and images of $B-b$ remain unchanged we have $\overline{HP}_*\cong HP_*$. This isomorphism is in fact an equality of homology groups. The only difference between $\overline{HP}$ and $HP$ is in the boundary maps for long exact sequences. The boundary maps in the long exact sequences in $\overline{HP}_*$ are defined using the boundary map in the complex and therefore are also multiplied by $\sqrt{2\pi i}$. We can now have a look at what happens to the transformation $\overline{\chi}_*: kk_*\to \overline{HP}_*$ defined in analogy to $\chi_*$. Since the maps $\alpha^n$ and $\beta^n$ used in Theorem \ref{Chern} for the definition of $\chi_*$ are defined using the boundary maps in the exact sequences in $HP_*$ for the universal and cone extensions, to obtain their analogues $\overline{\alpha}^n$, $\overline{\beta}^n$ in $\overline{HP}_*$ we have to multiply $\alpha^n, \beta^n$ by
$(2\pi i)^{n/2}$ and $(2\pi i)^{-n/2}$, respectively. Thus, if we define $\overline{\chi}_*:kk_*\to \overline{HP}_*$ in analogy to \eqref{chi2} using $\overline{\alpha}^n$, $\overline{\beta}^n$, then $\overline{\chi}$ has all the required properties - it is multiplicative, compatible with boundary maps and descends to a multiplicative transformation from $kk_*$, made $\Zz/2$-periodic, to $\overline{HP}_*$. In fact, under the natural identification $\overline{HP}_*=HP_*$ it is equal to the character $ch_*$ constructed above.\eremark

\end{document}